\magnification\magstep1 
\tolerance=10000

\def \nm#1{\left\|#1\right\|}

\def \sm#1{\sum_1^#1}
\def\expect{\hbox{{\rm I}\hskip -2pt {\rm E}}}
\def\superscript#1{^{\raise3pt\hbox{$\scriptstyle #1$}}}
\def \sstwo#1{_{\lower2pt\hbox{${\scriptstyle #1}$}}}

\def\Rn#1{\hbox{{\it I\kern -0.25emR}$\sp{\,{#1}}$}}
\def \seq#1#2{#1_1,\dots,#1_#2}
\def \qed {\vrule height6pt  width6pt depth0pt}

\def \ov{\overline}
\def \sup {{\rm sup\,}}
\def \span{{\rm span\,}}
\def \e{\epsilon }
\def \a{\alpha }
\def \b{\beta }
\def \s{\sigma}

\def \rarrow{\rightarrow}
\settabs 10 \columns
\def\ms{\medskip}
\def\bs{\bigskip}
\def\cen{\centerline}
\def\sm{\smallskip}
\def\ss#1{_{\lower3pt\hbox{$\scriptstyle #1$}}}

\def\o{\over}

\def \tr{{\rm trace\,}}

\def\Mbar{\ov{M}}
\def\Nbar{\ov{N}}

\def\Xt{\tilde{X}}

\def \E{{\Bbb E}}

\def \D{\Delta}

\def \nm#1{\left\|#1\right\|}

\def \seq#1#2{#1_1,\dots,#1_#2}

\def \sm#1#2{\sum_{#1=1}^#2}

\def \e{\epsilon}

\def \g{\gamma}
\def \d{\delta}
\def \t{\tau}


\catcode`\@=11
\font\tenmsy=msbm10
\font\sevenmsy=msbm7
\font\fivemsy=msbm5
\newfam\msyfam
\textfont\msyfam=\tenmsy  \scriptfont\msyfam=\sevenmsy
  \scriptscriptfont\msyfam=\fivemsy

\def\hexnumber@#1{\ifcase#1 0\or1\or2\or3\or4\or5\or6\or7\or8\or9\or
	A\or B\or C\or D\or E\or F\fi }

\edef\msy@{\hexnumber@\msyfam}

\def\Bbb{\ifmmode\let\next\Bbb@\else
 \def\next{\errmessage{Use \string\Bbb\space only in math mode}}\fi\next}
\def\Bbb@#1{{\Bbb@@{#1}}}
\def\Bbb@@#1{\fam\msyfam#1}

\catcode`\@=12

\cen{\bf Computing $p$--summing norms with few vectors}
\bs
\cen{ by William B.~Johnson\footnote*{Supported in part by NSF
\#DMS90-03550 and the U.S.-Israel Binational Science Foundation} and
Gideon Schechtman\footnote{**}{Supported in part by the U.S.-Israel
Binational Science Foundation}}

\

\noindent {\bf Abstract:} It is shown that the $p$-summing norm of any operator
with $n$-dimensional domain can be well-aproximated using only ``few" vectors in
the definition of the $p$-summing norm. Except for constants independent of
$n$  and $\log n$ factors, ``few" means $n$ if $1<p<2$ and $n^{p/2}$ if
$2<p<\infty$.

\

\

{\bf I. Introduction}

\

A useful result of Tomczak-Jaegermann [T-J, p. 143] states that the
2-summing norm of an operator $u$ of rank $n$ can be well-estimated by $n$
vectors; precisely (in the notation of [T-J, p. 140], which we follow
throughout),  $\pi_2(u) \le \sqrt{2} \pi_2^{(n)}(u)$.  No such result holds for
$\pi_1$; Figiel and Pelczynski [T-J, p. 184] showed that if $k_n$ satisfies
$\pi_1(u) \le C \pi_1^{(k_n)}(u)$ for all operators of rank $n$; $n=1,2,\dots$,
then $k_n$ grows exponentially in $n$.
The Tomczak result reduces
immediately to the case of operators whose domains are $\ell_2^n$.  Szarek [Sz]
proved that there is a $1$-summing analogue to this version of Tomczak's
theorem; namely, that $\pi_1(u) \le C \pi_1^{(n \log n)}(u)$ whenever $u$ is an
operator whose domain has dimension $n$.

In this paper we consider the case of $p$-summing operators.  In section III we
extend Szarek's result to the range $1<p<2$ (except that the power of $\log n$
is ``$3$" instead of ``$1$").  For $2<p<\infty$ we show that, up to powers of
$\log n$, $n^{p\o 2}$ vectors suffice to well-estimate the $p$-summing norm of
an operator from an $n$-dimensional space.  The power of $n$ is optimal, but
we do not know whether a $\log n$ term is needed in either result. These
results, as well as those in section IV, shed some light
on problems 24.10, 24.11, and 24.6 in [T-J].

In Section IV we  show that when $1 < p\neq 2 < \infty$, if $k_n$ satisfies
$\pi_p(u) \le C \pi_p^{(k_n)}(u)$ for all operators of rank $n$; $n=1,2,\dots$,
then $k_n$ grows faster than any power of $n$. .

Just as for Szarek, our main tools are sophisticated versions of embedding
$n$-dimensional subspaces of $L_p$ into $\ell_p^k$ with $k$ not too large.
While most of this background is at least implicit in [BLM] and [T], we need
more precise versions of such results than are stated in the current
literature.  The necessary material is developed in Section II.

Here we treat only the case of $p$-summing operators.  There is also an
extensive literature on related problems for $(p,q)$-summing operators; see
[N-T] for the older history and the recent papers  [DJ1], [DJ2], [J].   In
particular, Defant and Junge [DF2] show how results for $p$-summing operators
can be formally transformed into results for $(p,q)$-summing operators.

\

{\bf II. Preparations for the main result}
\ms

Before stating the basic entropy lemma for the main result, we set
some notation.

A {\it density} on a
probability space $(\Omega,\mu)$
is a strictly positive measurable function on $\Omega$ whose integral is
one. Given a set $A$, a metric $\d$ on $A$, and a positive number $t$,
$E(A,\d,t)$ is the minimal number of open balls of radius $t$ in the metric
$\d$ needed to cover $A$.

We also use notation (see,
for example, [T-J, p. 80]) commonly used in Banach space theory for measuring
the
expected value of the norm of Gaussian processes:   If $u : H\rarrow Z$ is a
linear
operator from a finite dimensional Hilbert space  $H$ into a normed space $Z$,
$\ell(u)^2$ is defined to be $\expect\Vert\sm{i}{m} g_i u(e_i)\Vert^2$, where
$\seq{e}{m}$ is any orthonormal basis for $H$ and $\seq{g}{m}$ are independent
standard Gaussian variables. $\ell$ is an ideal norm in the sense that if
$H^\prime$ is another finite dimensional Hilbert space, $Z^\prime$ is another
normed
space,  $T : H^\prime\rarrow Z$ and $S : Z\rarrow Z^\prime$ are linear
operators, then
$\ell(SuT) \le \Vert S\Vert \ell(u)\nm{T}$.  Suppose now that $\nu$ is a
probability measure on a finite set $A$ and $W$ is an $n$-dimensional subspace
of
the set of scalar valued functions on $A$.  Let $W_p$ denote $W$ under the
$L_p(\nu)$--norm and let $i_{p,r}^W$ be the formal identity mapping from $W_p$
onto
$W_r$ (when $r=2$ we   abuse notation by also regarding the
operator into $L_2(\nu)$).  Sudakov's lemma [Su], stated as Proposition 4.1
in
[BLM], gives the  entropy estimate
$$
\log E(B(W_p), \nm{\cdot}\sstwo{L_2(\nu)}, t) \le
C  \left({{\ell({i_{p,2}^W}^*)}\o t}\right)^2.
$$
The Pajor-Tomczak lemma [PT-J], stated as Proposition 4.2 in
[BLM], gives the  entropy estimate
$$
\log E(B(W_2), \nm{\cdot}\sstwo{L_p(\nu)}, t) \le
C  \left({{\ell({i_{2,p}^W})}\o t}\right)^2.
$$

The ideal properties of $\ell$ imply that if
$\nu$ (respectively, $\mu$) is a probability measure on the
finite set $A$ (respectively, $B$), and $Y$ (respectively $W$) is a space of
scalar functions on $A $ (respectively $B$), and $v : Y\rarrow W$ is a linear
operator which is an isometry from $Y_p$ onto $W_p$ and has norm at most $C$ as
an
operator from $Y_2$ into $W_2$, then $\ell({i_{p,2}^W}^*) \le
C\ell({i_{p,2}^Y}^*)$.

The entropy lemma we use is a variation on Propositions 4.6 and
7.2 in [BLM]. The result we need later is different from that in [BLM] since
we cannot replace the subspace $X$ of $L_p(\mu)$ by an (isomorphic or even
isometric) copy of $X$ in $L_p(\nu)$ but rather must move all of  $L_p(\mu)$
isometrically onto  $L_p(\nu)$.  Moreover, formally speaking, Proposition 7.2 is
only partly proved in [BLM] and contains some unclear statements (e.g., the
claim in the sentence immediately following (7.11) seems formally wrong and
should be adjusted slightly).  The accumulation of the adjustments needed to
obtain Proposition 2.1 below from the arguments in [BLM] required some effort on
our part, so we judged it worthwhile to outline proofs of the entropy
estimates we need.

\

 \proclaim Proposition 2.1. {Let $X$ be an $n$-dimensional subspace of
$L_p(\ov{N},\mu)$ for some probability measure $\mu$ on $\ov{N}=\{1,\dots,N\}$.
Then there is a density $\a$ on $(\ov{N},\mu)$ satisfying the
following: Put $\tilde{X} = \{x/{\a^{1\over p}} : \; x \in X \;\}$ and
let $B(\tilde{X}_r)$ be the closed unit ball of $\Xt$ in
$L_r(\ov{N},\a\,d\mu)$. Then for some constant $C$,
\item {(i)}
${\log E(B(\tilde{X}_p),\nm{\cdot}_{\infty},t) \le {C
{(p-1)\superscript{ \scriptscriptstyle{{p-2}\over 2}}}} n(\log
n)\superscript{1-{p\over 2}}(\log N)\superscript{\scriptscriptstyle{p\over
2}}t^{-p}}$,  for $1<p<2$. \vskip0.3cm
\item {(ii)}
$\nm{f}_{L_\infty} \le (2n)^{1/p} \nm{f}_{L_p(\a d\mu)}$, for
$1<p<2$.
\vskip0.3cm
\item {(iii)}
$\log E(B(\tilde{X}_p),\nm{\cdot}_{\infty},t) \le C(\log N)nt^{-2}$, for
$2\le p<\infty$.
\vskip0.3cm
\item {(iv)}
$\nm{f}_{L_\infty} \le (2n)^{1/2} \nm{f}_{L_p(\a d\mu)}$, for
$2<p<\infty$. }

\

\noindent{\bf Proof}. It is easy to reduce to the case of measures which
are strictly positive (i.e., for which all points of $\ov{N}$ have positive
$\mu$ measure).  The conclusion is invariant under change of density of the
original measure, so we can assume, without loss of generality, that $\mu$ is
the
uniform measure  on $\ov{N}$ (this simplifies slightly the notation below).

Lewis [L]
showed that there is a density $\b$ on $\ov{N}$ and an orthonormal (in
$L_2(\b d\mu)$) basis $\seq{f}{n}$ for
$Y = \{x/\a\sstwo{1}^{1/p}: \, x\in X \,\}$ so that $\sm{i}{n}
f_i^2 = n$.

The density $\a$ is ${{\b+1}\o 2}$.  Then
$\tilde{X}$ consists of all vectors of the form  $v(y)$ with $y$ in $Y$, where
$ v(y) = \left({\b \o \a}\right)^{1/p} y $. The linear operator  $v$  defines
an
isometry from $L_p(\b d\mu)$ onto $L_p(\a d\mu)$ and  has norm at
most $2^{1/p}$ as an operator from $L_2(\b d\mu)$ into $L_2(\a d\mu)$.    As
mentioned before the statement of Proposition 2.1, Sudakov's lemma gives the
entropy estimate $$
\log E(B(\tilde{X}_p), \nm{\cdot}\sstwo{L_2(\a\mu)}, t) \le C t^{-2} p n
K(X)^2.
$$

Using the fact that
$\seq{({\b\o\a})^{1/2}f}{n}$ is orthonormal in $L_2(\a d\mu)$
and the Maurey-Khintchine inequality, we get for all $1 \le q < \infty$:
 $$
\eqalign{\ell^2({i_{2,q}^{{\tilde{X}}}}) & =
\expect\nm{\sm{i}{n} g_i \left({\b\o \a}\right)^{1/2}f_i }_{L_q(\a d\mu)}^2
\cr &
\le C q \nm{(\sm{i}{n}  \left({\b\o \a}\right)f_i^2)^{1/2} }_{L_q(\a d\mu)}^2
\cr &
\le 2 C q n.}
$$

As mentioned before the statement of Proposition 2.1, the Pajor-Tomczak
lemma gives the entropy estimate
$$
\log E(B(\tilde{X}_2), \nm{\cdot}\sstwo{L_q(\a d\mu)}, t) \le C  t^{-2} q
n.\eqno{(+)} $$

Pick $q=\log 2N$; then, since $\a > 1/2$,
$\nm{\cdot}_\infty \le e \nm{\cdot}_{L_q(\a d\mu)}$.  Since $B(\tilde{X}_p)
\subset B(\tilde{X}_2)$ for $p \ge 2$,
this gives (iii). The Lewis change of density forces, for $f$ in $Y$,
$\nm{f}_{L_\infty} \le n^{1/2} \nm{f}_{L_p(\b d\mu)}$  (see e.g.
Lemma 7.1 in [BLM]).  Since
${\b \o \a} \le 2$, we have for  $f$ in $\tilde{X}$ that
$\nm{f}_{L_\infty} \le (2n)^{1/2} \nm{f}_{L_p(\a d\mu)}$.  This gives (iv).

To deal with the case  $1<p<2$, we refer to the proof of Proposition 7.2 (ii)
in
[BLM]. By applying H\"older's inequality and a clever duality argument, one
obtains formally from (+),  for
$1 \le t \le 2n$, that
$$
\log E(B(\tilde{X}_p), \nm{\cdot}\sstwo{L_2(\a d\mu)}, t)
\le C (p-1)^{-1} ({C\o t})^{2p/(2-p)}  n\log n.\eqno{(++)}
$$
 Using, for $1< s < t$, the obvious inequality
$$\log E(B(\tilde{X}_p), \nm{\cdot}\sstwo{L_\infty},t)
\le \log E(B(\tilde{X}_p), \nm{\cdot}\sstwo{L_2(\a d\mu)}, s) +
\log E(B(\tilde{X}_2), \nm{\cdot}\sstwo{L_\infty}, t/s),
$$
 (++), and (iii) in the statement of the proposition, we obtain (i) for
$1\le t \le 2n$ by minimizing over
$s$.  Now the Lewis change of density forces, for $f$ in $Y$,
$\nm{f}_{L_\infty} \le n^{1/p} \nm{f}_{L_p(\b d\mu)}$.  Since
${\b \o \a} \le 2$, we have for  $f$ in $\tilde{X}$ that
$\nm{f}_{L_\infty} \le (2n)^{1/p} \nm{f}_{L_p(\a d\mu)}$.  This gives (ii) as
well as  (i) when $t>2n$.  \hfill\qed

\
The following proposition and its proof is an adjustment of results from
Talagrand's paper [T]. (The idea of ``splitting the large atoms", used also in
 [T], is
due to the authors.)
\
\proclaim Proposition 2.2. {Let $X$ be an $n$-dimensional subspace of
$L_p(\ov{N},\t)$ for some probability measure $\t$ on $\ov{N}=\{1,\dots,N\}$.
Then there are $N\le M\le {3\over 2} N$ and a probability measure $\nu$ on
$\ov{M}=\{1,\dots,M\}$ so that:
\vskip0.1cm
\item{(i)}
There is a partition $\{\s_1,\dots,\s_n\}$ of $\ov{M}$ with ${{
\sum_{i\in \s_j} \nu\{i\} =\t\{j\}}}$ for $j=1,\dots,n$.
\item {(ii)}
${{\E \, \sup \{ \bigl| \sum_{i=1}^M g_i\nu\{i\}|y_i|^p\, \bigr| : y\in Y,
||y|| \le 1 \}\le
C
{(p-1)\superscript{ \scriptscriptstyle{{p-2}\over 4}}} ({n\o
N})\superscript{{1\o 2}}(\log n)
\superscript{{{6-p}\over
4}}(\log N)\superscript{\scriptscriptstyle{p\over 4}}}}$, for $1<p<2$,
where $y_1,\dots,y_n$ are the coordinates of the vector $y$
and $Y$ is the image of $X$ under the natural isometry $J_p$ from
$L_p(\ov{N},\mu)$
into $L_p(\ov{M},\nu)$, defined by $(J_px)_i = x_j$ if $i\in \s_j$.
\vskip0.1cm
\item {(iii)}
${{\E \, \sup \{ \bigl| \sum_{i=1}^M g_i\nu\{i\}|y_i|^p\, \bigr| : y\in Y,
||y|| \le 1 \}\le
C_pn^{p/4}N^{-1/2}\log n (\log N)^{1/2}}}$,\hfill\break
for $2<p<\infty$. $C_p$ can be taken to be $Cp^2 2^{p/2}$.
\vskip0.1cm}

\noindent{\bf Proof}.  It is easy to reduce to the case of measures which
are strictly positive. Next, note that if the proposition  is true for  {\sl
one}
strictly positive probability measure on $\ov{N}$, then it is true for all of
them. This is because the left hand side of (ii) is invariant under a change of
density $\phi$ if we replace the  subspace $Y$ of $L_p(\nu)$ with its image
under
the natural isometry from $L_p(\nu)$ onto $L_p(\phi\,d\nu)$, defined by $Tf =
f/\phi^{1/p}$.  Thus we can assume that  $\t$ is the measure $\a d\mu$ given by
the conclusion of Proposition 2.1.

Splitting the atoms of $\t$
of mass larger than $4/N$ into pieces each of size between $2/N$ and $4/N$
produces $\ov{M}$, the measure $\nu$, and, {\it a fortiori}, the space $Y$
along
with the isometry $J=J_p$; (i) is thus satisfied.  Since $J$ also defines an
isometry $J_r$ from $L_r(\ov{N},\t)$ into $L_r(\ov{M},\nu)$ for all $0 < r \le
\infty$,
the conclusion of Proposition 2.1 remains true for the measure space
$(\ov{M},\nu)$ (where of course $\tilde{X}$ is replaced by $Y$).

 Let $\d$ be the natural distance
associated with the Gaussian process appearing in (ii), defined for $y, z$ in
$Y$ by
$$
\d(y,z) = \left(\sm{i}{M} [\t\{i\}(|y_i|^p - |z_i|^p)]^2\right)^{1/2}.
$$
Let $1 < p < 2$, fix $y,z$ in  $B(Y_p)$, and set $u_i = |y_i|\vee |z_i|$. Then
$$
\eqalign{\d(y,z)^2 & \le \sm{i}{M} \nu\{i\}^2 p^2 u_i^{2p-2} |y_i - z_i|^2
\cr
& \le \nm{y-z}_\infty^p 4p^2 N^{-1} \sm{i}{M} \nu\{i\}  u_i^{2p-2} |y_i -
z_i|^{2-p}
\cr
& \le  4p^2 N^{-1} \nm{y-z}_\infty^p
(\sm{i}{M} \nu\{i\}  u_i^p)^{2(p-1)/p} (\sm{i}{M} \nu\{i\} |y_i -
z_i|^p)^{(2-p)/p}
\cr
& \le 2^6  N^{-1} \nm{y-z}_\infty^p.}
$$
Thus by Proposition 2.1 (ii) we get that the $\d$-diameter of $B(Y_p)$ is less
than
$2^4 n^{1/2}N^{-1/2}$ and from Proposition 2.1 (i) that:
$$\eqalign{
\log E(B(Y_p),\d,t) & \le \log E(B(Y_p),\nm{\cdot}_\infty^{p/2}, 2^{-3} N^{1/2}
t)
\cr &
\le \log E(B(Y_p),\nm{\cdot}_\infty, 2^{-6/p} N^{1/p} t^{2/p})
\cr &
\le C
{(p-1)\superscript{ \scriptscriptstyle{{p-2}\over 2}}} n(\log
n)\superscript{1-{p\over
2}}(\log N)\superscript{\scriptscriptstyle{p\over 2}} N^{-1} t^{-2}.}
$$
The last inequality in this last display requires $t\ge 2^3 N^{-1/2}$;
for $0 < t < 2^3 N^{-1/2}$ use volume considerations in the $n$-dimensional
space
$B(Y_\infty)$ to get
$$\eqalign{
\log E(B(Y_p),\d,t)  &
\le \log E(B(Y_p),\nm{\cdot}_\infty,1) +
\log E(B(Y_\infty),\nm{\cdot}_\infty, 2^{-6/p} N^{1/p} t^{2/p})
\cr &
\le  C
{(p-1)\superscript{ \scriptscriptstyle{{p-2}\over 2}}} n(\log
n)\superscript{1-{p\over
2}}(\log N)\superscript{\scriptscriptstyle{p\over 2}} +
C n \log (C N^{-1} t^{-2}).    }
$$
By Dudley's theorem (see, e.g., [MP, p. 25]),
$$\eqalign{
\E \; \sup &\{ \bigl| \sum_{i=1}^M g_i\nu\{i\}|y_i|^p\, \bigr| : y\in Y, \;\;\;
||y|| \le 1 \}
 \cr &
\le
2^4 n^{1/2}N^{-1/2} +
C
{(p-1)\superscript{ \scriptscriptstyle{{p-2}\over 4}}} n^{1/2}(\log
n)\superscript{{{2-p}\over
4}}(\log N)\superscript{\scriptscriptstyle{p\over 4}}N^{-1/2}
\cr
&
\quad +
C n^{1/2} \int_0^{2^3N^{-1/2}} \log^{1/2} (CN^{-1}t^{-2}) \, dt
\cr &
\quad +
C
{(p-1)\superscript{ \scriptscriptstyle{{p-2}\over 4}}} n^{1/2}(\log
n)\superscript{{{2-p}\over
4}}(\log N)\superscript{\scriptscriptstyle{p\over 4}}N^{-1/2}
\int_{2^3N^{-1/2}}^{2^4 n^{1/2}N^{-1/2}} t^{-1} \, dt
\cr
& \le
C
{(p-1)\superscript{ \scriptscriptstyle{{p-2}\over 4}}} n^{1/2}(\log
n)\superscript{{{6-p}\over
4}}(\log N)\superscript{\scriptscriptstyle{p\over 4}}N^{-1/2}.
}
$$
This proves (ii).

To prove (iii), assume now $2<p<\infty$. Fix $y,z$ in  $B(Y_p)$, and set $u_i =
|y_i|\vee |z_i|$. Then
$$
\eqalign{\d(y,z)^2 & \le \sm{i}{M} \nu\{i\}^2 p^2 u_i^{2p-2} |y_i - z_i|^2
\cr
& \le \nm{y-z}_\infty^2 4p^2 N^{-1} \sm{i}{M} \nu\{i\}  u_i^{2p-2}
\cr
& \le  4p^2 N^{-1} \nm{y-z}_\infty^2 \nm u_\infty^{p-2}
\sm{i}{M} \nu\{i\}  u_i^p
\cr
& \le 4p^2 2^{p/2}{{n^{(p-2)/2}}\o N}\nm{y-z}_\infty^2,}
$$
where the last inequality follows from Proposition 2.1 (iv).
Thus the $\d$ diameter
of $B(Y_p)$ is less than
$4p2^{p/4} n^{p/4}N^{-1/2}$ and Proposition 2.1 (iii) implies:
$$\eqalign{
\log E(B(Y_p),\d,t) & \le \log E(B(Y_p),\nm{\cdot}_\infty,
p^{-1}2^{-(p+4)/4}n^{-(p-2)/4}
 N^{1/2} t)
\cr &
\le C
p^2 2^{p/2}n^{p/2}N^{-1}(\log N)t^{-2}}
$$
as long as $t\ge  p2^{(p+4)/4}n^{(p-2)/4}
 N^{-1/2}$. For smaller $t$ we get by the usual volume considerations,
$$
\eqalign{
\log E(B&(Y_p),\d,t)  \cr
&
\le \log E(B(Y_p),\nm{\cdot}_\infty,1) +
\log E(B(Y_\infty),\nm{\cdot}_\infty, p^{-1}2^{-(p+4)/4}n^{-(p-2)/4}
 N^{1/2} t)
\cr &
\le  Cp^2 2^{p/2}n^{p/2}N^{-1}(\log N) +
C n \log (C p 2^{p/4}n^{(p-2)/4}N^{-1/2} t^{-1}).    }
$$
By Dudley's theorem,
$$\eqalign{
\E \; \sup &\{ \bigl| \sum_{i=1}^M g_i\nu\{i\}|y_i|^p\, \bigr| : y\in Y, \;\;\;
||y|| \le 1 \}
\cr &
\le
Cp2^{p/4} n^{p/4}N^{-1/2} + Cp^22^{p/2}n^{(p-1)/2}N^{-1}(\log N)^{1/2}\cr
&\quad +
C n^{1/2} \int_0^{p2^{(p+4)/4}n^{(p-2)/4}
N^{-1/2}} \log^{1/2} (C p 2^{p/4}n^{(p-2)/4}N^{-1/2} t^{-1})\, dt\cr
&  \quad +
Cp 2^{p/4}n^{p/4}N^{-1/2}(\log N)^{1/2}\int_{p2^{(p+4)/4}n^{(p-2)/4}N^{-1/2}}^
{4p2^{p/4} n^{p/4}N^{-1/2}}t^{-1}\, dt.
\cr}$$
For a fixed $p$ the last term is dominating and one gets
$$
\E \; \sup \{ \bigl| \sum_{i=1}^M g_i\nu\{i\}|y_i|^p\, \bigr| : y\in Y, \;\;\;
||y|| \le 1 \}
\le C_pn^{p/4}N^{-1/2}\log n (\log N)^{1/2}
$$
where $C_p$ can be taken to be $Cp^2 2^{p/2}$.
\hfill\qed

\

\proclaim Corollary 2.3. {Let $X$ be an $n$-dimensional subspace of
$L_p(\ov{N},\t)$ for some probability measure $\t$ on $\ov{N}=\{1,\dots,N\}$
and let $L_p(\Mbar,\nu)$,  $J$, and $Y$ be given from Proposition 2.2.  Then
there  is a partition $M_1\cup M_2$  of
$\Mbar$ into two  sets of cardinality at most ${7\o 8}N$ such that
for each $y$
in $Y$ and $j=1,2$:
{\item{(i)} $\nm{1_{M_j}y}\ss{L_p(\Mbar,\nu)}^p \le
\left({1/2}+ C
{(p-1)\superscript{ \scriptscriptstyle{{p-2}\over 4}}}
({n\o N})\superscript{{1\o 2}}(\log n)\superscript{{{6-p}\over
4}}(\log N)\superscript{\scriptscriptstyle{p\over 4}}\right)
\nm{y}\ss{L_p(\Mbar,\nu)}^p,$} \break
\vskip0.2cm
\noindent
when $1<p<2$; while
{\item{(ii)} $\nm{1_{M_j}y}\ss{L_p(\Mbar,\nu)}^p \le
\left({1/2}+ C_p
\left({n\superscript{{p\o 2}}\o N}\right)\superscript{{1\o 2}}\log
n(\log N)\superscript{\scriptscriptstyle{1\over
2}}\right) \nm{y}\ss{L_p(\Mbar,\nu)}^p,
$ for $2<p<\infty.$}
\break
Moreover, (i) and (ii) hold for most such partitions of $\Mbar$.
}

\

\noindent{\bf Proof}.  First, notice that (ii) in Proposition 2.2 still holds
if we substitute independent Rademacher functions for the Gaussian variables
$g_i$ (and replace $C$ by, e.g., $\sqrt{\pi\o 2} C$).  This follows from a
standard contraction principle.  Consequently, if we again enlarge $C$,
$$
{{ \sup \{ \bigl| \sum_{i=1}^M \e_i\nu\{i\}|y_i|^p\, \bigr| : y\in Y,
||y|| \le 1 \}\le
C
{(p-1)\superscript{ \scriptscriptstyle{{p-2}\over 4}}} ({n\o
N})\superscript{{1\o
2}}(\log n)\superscript{{{6-p}\over 4}}(\log
N)\superscript{\scriptscriptstyle{p\over 4}}}} $$
holds for most choices of signs $\e_i=\pm 1$. Since also for most choices of
signs the  difference between the number of plus signs and minus signs  is less
than $M/8$, (i) follows. (ii) follows similarly. \hfill\qed

\

{\bf III. Computing $p$-summing norms}
\ms
Given a linear operator $u : X\to Y$ of finite rank, $1\le q\le \infty$, and
 positive integers $n$,
 $k$, define
$$\nu_q^{(n,k)}(u) = \inf\biggl{\{}\sm{i}{k} \nu_q^{(n)}(u_i) : \; u= \sm{i}{k}
u_i \biggr{\}},
$$
where
$$
\nu_q^{(n)}(v) = \inf \left\{ \nm{A}\nm{w}\nm{B} ; \; A :  X\to \ell_\infty^n
;  w : \ell_\infty^n \to \ell_q^n \hbox{ diagonal, } B : \ell_q^n \to Y,
  v = BwA \right\}.
$$

In Tomczak's terminology [T-J, p. 181], $\nu_q^{(n,1)}= \nu_q^{(n)}$, while
$\lim_{k\to \infty} \nu_q^{(n,k)} = \hat{\nu}_q^{(n)}$ gives the cogradation
which is dual to the natural gradation $\pi_p^{(n)}$ of the $p$-summing
norm  [N-T, Theorem 24.2] (or something like that!).

\

\proclaim Proposition 3.1. {Let $n\le N$ be positive integers; $u : X\to Y$ a
linear operator with $X$ finite dimensional and $\dim (Y) \le n$.  Then,
putting
$q=p/(p-1)$,
 \item{(i)}
For $1<p<2$ ,
$$
\nu_q^{({7\o 8}N,2)}(u) \le \left(1+C
{(p-1)\superscript{ \scriptscriptstyle{{p-2}\over 4}}}
\left({n\o N}\right)\superscript{{1\o 2}}(\log n)\superscript{{{6-p}\over
4}}(\log N)\superscript{\scriptscriptstyle{p\over 4}}\right)\nu_q^{(N,1)}(u).
$$
\item{(ii)}
For $2<p<\infty$,
$$\nu_q^{({7\o 8}N,2)}(u) \le \left(1+C_p
\left({n\superscript{{p\o 2}}\o N}\right)\superscript{{1\o 2}}\log
n(\log N)\superscript{\scriptscriptstyle{1\over
2}}\right)\nu_q^{(N,1)}(u).
$$
}

\

\noindent{\bf Proof}.  For some probability measure $\t$ on $\Nbar$, we can
take $A :Y^*\to L_p(\Nbar,\t)$, \break $B : L_1(\Nbar,\t) \to X^*$, so that
$\nm{A} \nm{B} =  \nu_q^{(N)}(u)$ and
$u^* = Bi_{p,1}A$. Apply Proposition 2.2 to
the subspace $AY$ of $L_p(\Nbar,\t)$ to  get the measure space
$L_p(\Mbar,\nu)$ and the natural isometric embedding
$J_p : L_p(\Nbar,\t) \to L_p(\Mbar,\nu)$. By
Corollary 2.3,  we get a partition $M_1\cup M_2$  of
$\Mbar$ into two  sets of cardinality at most ${7\o 8}N$ such that for each $y$
in $Y$, $j=1,2$, and in the case $1<p<2$:
 $$\nm{1_{M_j}J_pAy}\ss{L_p(\Mbar,\nu)}^p \le
\left({1/2}+ C
{(p-1)\superscript{ \scriptscriptstyle{{p-2}\over 4}}}
({n\o N})\superscript{{1\o 2}}(\log n)\superscript{{{6-p}\over
4}}(\log N)\superscript{\scriptscriptstyle{p\over 4}}\right)
\nm{Ay}\ss{L_p(\Nbar,\t)}^p.
$$
 Denote for $j=1,2$ the injection from $L_p(M_j,\nu_{|M_j})$ to
$L_1(M_j,\nu_{|M_j})$ by $i_{p,1}^j$ and let $P$ be the conditional
expectation projection from $L_1(\Mbar,\nu)$ onto $J_1[L_1(\Nbar,\t)]$
followed by $J_1^{-1}$. Thus
$u^* = BPi_{p,1}^11_{M_1}J_1A + BPi_{p,1}^21_{M_2}J_1A$
and
$$\eqalign{\nu_q^{({7\o 8}N,2)}(u) &\le
\sm{j}{2} \nu_q^{({7\o 8}N)}([BPi_{p,1}^j1_{M_j}J_1A]^*)\cr
&\le \sm{j}{2}\nm{1_{M_j}J_1A}\bigl\Vert i_{p,1}^j\bigr\Vert \bigl\Vert
BP\bigr\Vert\cr
&\le \left({{1\o 2}}+ C
{(p-1)\superscript{ \scriptscriptstyle{{p-2}\over 4}}}
({n\o N})^{1\o 2}(\log n)\superscript{{{6-p}\over
4}}(\log N)\superscript{\scriptscriptstyle{p\over 4}}\right)\superscript
{{1\o p}} \nm{A}\nm{B}\sm{j}{2}\nu(M_j)^{1\o q} \cr
&\le \nm{A}\nm{B}
\left(1+2C
{(p-1)\superscript{ \scriptscriptstyle{{p-2}\over 4}}}
({n\o N})\superscript{{1\o 2}}(\log n)\superscript{{{6-p}\over
4}}(\log N)\superscript{\scriptscriptstyle{p\over 4}}\right)\superscript{1\o
p}\cr
&\le \nm{A}\nm{B}
\left(1+2C
{(p-1)\superscript{ \scriptscriptstyle{{p-2}\over 4}}}
({n\o N})\superscript{{1\o 2}}(\log n)\superscript{{{6-p}\over
4}}(\log N)\superscript{\scriptscriptstyle{p\over 4}}\right).}
$$
This completes the proof when $1<p<2$; the other case
is similar.
\hfill\qed

\

\proclaim Theorem 3.2. {Suppose that $\dim (X) \le n$, $u : X\to Y$ is a
linear operator and $\e > 0$.  Then,
$$
\pi_p(u) \le (1 + \e) \pi_p^{(m)}(u),
$$
as long as
\item{(i)}
$1<p<2$ and \space\space$m \ge {K  {(p-1)\superscript{
\scriptscriptstyle{{p-2}\over
2}} \e^{-2}}} n (\log n)\superscript{ \scriptscriptstyle{{6-p}\over
2}}\left(\log \left(
{(p-1)\superscript{ \scriptscriptstyle{{p-2}\over 2}} \e^{-2}}
n\right)\right)\superscript{p \o 2} $ \space\space
for some absolute constant $K$,\hfill
\break
\noindent or
\vskip0.1cm
\item{(ii)}
$2<p<\infty$ and \space\space$m \ge K_p\e^{-2}n\superscript{
\scriptscriptstyle{ p\over 2}} (\log n)^2\log(\e^{-2}n\superscript{
\scriptscriptstyle {p\over 2}})$.
}

\

\noindent{\bf Proof}.
Without loss of generality, we can assume that $\dim(Y) \le n$. By duality
[T-J,
Theorem 24.2], it is enough to prove that
$$
\hat\nu_q^{(m)}(v) \le (1 + \e)
\nu_q^N(v)
$$
 for all $v : Y \to X$ and all positive integers $N\ge n$. Iterating
Proposition 3.1, we get for all $k$ (with $({7\o 8})^kN \ge n$) and for
$1<p<2$ that
$$\eqalign{
\nu_q^{([{7\o 8}]^kN,2^k)}&(u) \le \cr
&
\prod_{j=1}^k
\left(1+C
{(p-1)\superscript{ \scriptscriptstyle{{p-2}\over 4}}}
\left({n\o ({7\o 8})^{j-1}N}\right)\superscript{{1\o 2}}(\log
n)\superscript{{{6-p}\over
4}}(\log \left[ \left(\scriptstyle{7\o
8}\right)^{j-1}N\right])\superscript{\scriptscriptstyle{p\over 4}}
\right)\nu_q^{(N)}(u).} $$
The product on the right hand side of the above inequality is smaller than
$1+\e$ as long as
$$
{(p-1)\superscript{ \scriptscriptstyle{{p-2}\over 4}}}
\left({n\o ({7\o 8})^{k}N}\right)\superscript{{1\o 2}}(\log n)
\superscript{{{6-p}\over
4}}(\log \left[ \left(\scriptstyle{7\o 8}\right)^{k}N\right])
\superscript{\scriptscriptstyle{p\over
4}} \le \d \e,
$$
(where $\d = \d(C)$ is an appropriate positive constant). Put
$m=\left(\scriptstyle{7\o 8}\right)^{k}N$; then, as long as
 \space\space $m \ge {\d'  {(p-1)\superscript{ \scriptscriptstyle{{p-2}\over
2}}
\e^{-2}}} n (\log n)\superscript{ \scriptscriptstyle{{6-p}\over 2}}\left(\log
\left(
{(p-1)\superscript{ \scriptscriptstyle{{p-2}\over 2}} \e^{-2}}
n\right)\right)\superscript{p \o 2} $, \space\space
$$
\hat\nu_q^{(m)}(u) \le \nu_q^{(m,2^k)}(u) \le \nu_q^{(N)}(u).
$$
This completes the proof when $1 < p < 2$; the case $2 < p < \infty$ is
similar.
\hfill\qed

\

{\bf Remark.} As we have presented it, the proof of Theorem 3.2 does not
recapture the result of Szarek mentioned in the introduction.  Actually, our
approach does work when $p=1$ and the technical difficulties are easier in this
case because the entropy considerations of Section II are not needed.

\

\def\at{\tilde{\a}}
{\bf IV. Examples and concluding remarks}
\ms
For $p>2$, $\pi_p^{(k)}(\ell_2^n)\le k^{1\o p}$, while
$\pi_p(\ell_2^n)\ge \sqrt{n/p}$\ [N-T, Theorem 10.2]. Consequently, Theorem 3.2
is precise except for the $\log n$ terms.

It is natural to ask what value of $k$ is needed for $\pi_p^{(k)}(u)$ to
well-estimate $\pi_p(u)$ for a general operator $u$ of rank $n$.  When
$p=1$, Figiel-Pe\l czynski [T-J, p.184] checked that $k$ must be exponential in
$n$.  The authors and J. Bourgain checked that a result of Bourgain's [B]
yields that for $1 < p \neq2 <\infty$, $k$ grows faster than any power of $n$.

\proclaim Proposition 4.1. Let $1 < p \neq2 <\infty$ and $C<\infty$.  Suppose
that for each $s=1,2,\dots$,\hfill\break
 $k_s$ satisfies
$$
\pi_p(u) < C\pi_p^{(k_s)}(u)
$$
for all operators $u$ of rank at most $s$.  Then for all $K<\infty$,
$k_s s^{-K} \rarrow \infty$ as $N\rarrow \infty$.

\noindent{\bf Proof:} Fix $1 < p \neq2 <\infty$, $K$, $C$, and let $\d>0$ with
$\d K <1$.  Given $N=2^n$ for some $n$, we identify $L_p^N$ with $L_p(G)$,
where $G$ is the group $\{-1,1\}^n$ with normalized Haar measure, d$ g$.

Let $E = \span\{w_S : |S| \ge n-m \}$ where ${m\o n} \log {n\o m} \sim \d$;
so $\dim E \sim ({n\o m})^m < N^\d$.
Here we follow Bourgain's notation [B]; for $S\subset \{1,\dots,n\}$,
$w_S = \prod_{i\in S} r_i$, with $r_i$ the $i$-th coordinate projection
(Rademacher) on $G$. Let $j_{\infty,p}^E$ be the formal identity from
$E_{\infty}$ to $L_p^N$.  We shall use Bourgain's result [B] that if $T$ is an
operator on $L_p^N$ which is the identity on $E$ and $\nm{T} < C$, then $\tr{T}
\sim N$ (meaning $|\tr{(I-T)}| =o(N)$), to prove that if
$$
\hat{N}_p^{(k)}(j_{\infty,p}^E) < C\nu_p(j_{\infty,p}^E)\,\,\,\,(=C),
$$
then for large $N$, $\,\,\, k>N^{\d K} \ge (\dim E)^K$. This gives the dual
form of the conclusion of Proposition 4.1.

For notational convenience, set $\a=j_{\infty,p}^E$ and suppose that for
certain $k$ we have $\a=\sum_i \a_i$ with $\sum_i \nu_p^{k} (\a_i) < C$.
This means that there are factorizations
$$
E\;{{\t_i}\atop\raise9pt\hbox{$\longrightarrow$}}\;
\ell_\infty^k\;
{{\D_i}\atop\raise9pt\hbox{$\longrightarrow$}}\;\ell_p^k \;
{{\g\ss{i}}\atop\raise9pt\hbox{$\longrightarrow$}}\;L_p^N
$$
of $\a_i$ with $\nm{\t_i}=\nm{\D_i}=1$, $\D_i$ diagonal, and
$\sum\ss{i} \nm{\g\ss{i}}<C$.  This diagram also gives that
$\sum\ss{i}\nu_1(\a_i) < Ck^{1-{1\o p}}$.
Extend $\t_i$ to a map $\tilde{\t}_i : L_\infty^N \to
\ell_{\infty}^k$ with $\nm{\tilde{\t}_i} =1$, set
$\tilde{\a}_i = \g\ss{i} \D_i \tilde{\t}_i$, and let
$\tilde{\a}=\sum\ss{i} \at_i$.  Then
$$
\nu_1(\at) \le \sum_i \nu_1(\at_i) < C k^{1-{1\o p}}\;\;
\hbox{and} \;\;
\pi_p(\at) = \nu_p(\at) < C.
$$
Now replace $\at$ by its average $\b$ over the group
$G$, defined by
$$
\b = \int_G T_g\at T_g \hbox{d}\,g \quad (g^{-1} = g \;\;
\hbox{in} \;\; G).
$$

The operator $\b$ is translation invariant (a multiplier) and
satisfies the same conditions as $\at$; namely,
$$\b_{|E} = \a,\quad \nu_1(\b) < C k^{1-{1\o p}},
\quad \pi_p(\b) < C.
$$
Since $\b$ is translation invariant, Haar measure on $G$ is a
suitable Pietsch measure, which means that
$\nm{\b i_{p,\infty}}<C$.  Thus $\tr (\b i_{p,\infty}) \sim N$
by Bourgain's result [B].  However,
$$
|\tr (\b i_{p,\infty})| \le \nu_1(\b i_{p,\infty}) \nm{i_{p,\infty}}
< C k^{1-{1\o p}}N^{1\o p},
$$
which is $o(N)\;$ if $\;k\le N^{\d K}$. \hfill \qed

\vfill
\eject

\

\centerline{\bf References}

\

\item{[B]} J. Bourgain, {\sl A remark on the behaviour of $L^p$-multipliers and
the range of operators acting on $L^p$-spaces,} {\bf Israel J. Math. 79}
(1992), 1--11.

\item{[BLM]} J. Bourgain, J. Lindenstrauss and V.~D. Milman, {\sl
Approximation of zonoids by zonotopes,} {\bf Acta Math.162} (1989), 73--141.

\item{[DJ1]} M. Defant and M. Junge, {\sl Absolutely summing norms with $n$
vectors,}

\item{[DJ2]} M. Defant and M. Junge, {\sl On absolutely summing operators with
application to the
$(p,q)$-summing norm with few vectors,}

\item{[D]}  R.~M. Dudley, {\sl The sizes of compact subsets of Hilbert space
and continuity of Gaussian processes,} {\bf J. Functional Analysis 1}
(1967), 290--330.

\item{[J]} G.~J.~O. Jameson, {\sl The number of elements required to determine
$\pi_{2,1}$,}

\item{[L]} D.~R. Lewis, {\sl  Ellipsoids defined by Banach ideal norms,} {\bf
Mathematika 26} (1979), 18--29.

\item{[MP]}  M.~B. Marcus and G. Pisier, {\sl Random Fourier series with
applications to harmonic analysis,} {\bf Annals of Math. Studies 101}
(1981), Princeton University Press, N. J.

\item{[Sc]} G. Schechtman, {\sl More on embedding subspaces of  $L_p$
into  $l_r^n$,} {\bf Comp. Math. 61} (1987), 159--170.

\item{[Sz]} S.~J. Szarek, {\sl  Computing summing norms and type constants
on few vectors,} {\bf Studia Math. 98} (1990), 147-156.

\item{[T]} M. Talagrand, {\sl Embedding subspaces of  $L_1$ into
$l_1^N$,} {\bf  Proc. A.M.S. 108} (1990), 363--369.

\item{[T-J]} N. Tomczak-Jaegermann, {\sl Banach-Mazur distances and
finite-dimensional operator ideals,} Pitman Monographs and Surveys
in Pure and Applied Mathematics {\bf 38}, Longman 1989.

{\baselineskip=14pt

{\  }

{\ }

 {\obeylines \parindent=0pt \lineskip=2pt
Department of Mathematics,
Texas A\&M University,
College Station, TX 77843, USA
e-mail: wbj7835@venus.tamu.edu

\

Department of Theoretical Mathematics,
The Weizmann Institute of Science,
Rehovot, Israel
e-mail: mtschech@weizmann.weizmann.ac.il

{\ } }}

\bye